\documentclass[psamsfonts,reqno]{amsart}

\usepackage{amssymb}
\usepackage{amscd}
\usepackage{eucal}
\usepackage{psfig}
\usepackage[dvips]{graphicx}
\usepackage[small]{caption}

\newcommand{\beq}{\begin{equation}}
\newcommand{\eeq}{\end{equation}}

\newcommand{\sds}{\strut\displaystyle}

\newcommand{\N}{\mathbb{N}}

\newcommand{\cE}{\mathcal{E}}
\newcommand{\cF}{\mathcal{F}}

\newcommand{\cK}{\mathcal{K}}
\newcommand{\cI}{\mathcal{I}}

\newcommand{\cN}{\mathcal{N}}

\newcommand{\sfrac}[2]{{\vphantom1\smash{\lower.5ex\hbox{\small$#1$}}\over
        \vphantom1\smash{\raise.4ex\hbox{\small$#2$}}}} 

\newcommand{\barg}{{\bar g}}
\newcommand{\smCs}{{\mbox{\tiny Cons.}}}
\newcommand{\EA}{\mbox{E}}
\newcommand{\Tr}{\mbox{Tr}\,}
\newcommand{\Ra}{{\mbox{Range}\,}}

\def\staccrel#1#2{\mathrel{\mathop{#1}\limits_{#2}}}

\newtheorem{theorem}{Theorem}
\newtheorem*{theorem-nonumber}{Theorem}

\newtheorem{lemma}{Lemma}

\theoremstyle{remark}
\newtheorem*{remark}{Remark}

\theoremstyle{definition}
\newtheorem*{problem}{Problem}

\begin{document}

\title[Information associated with FIE's]{Metric and Probabilistic Information Associated with Fredholm Integral Equations of the First Kind}

\author[E. De Micheli]{Enrico ~De Micheli}
\address[E. De Micheli]{IBF -- Consiglio Nazionale delle Ricerche \\ Via De Marini, 6 - 16149 Genova, Italy}
\email[E.~De Micheli]{demicheli@ge.cnr.it}

\author[G. A. Viano]{Giovanni Alberto ~Viano}
\address[G. A. ~Viano]{Dipartimento di Fisica - Universit\`a di Genova,
Istituto Nazionale di Fisica Nucleare - sez. di Genova\\
Via Dodecaneso, 33 - 16146 Genova, Italy}
\email[G.A.~Viano]{viano@ge.infn.it}

\begin{abstract}
The problem of evaluating the information associated with Fredholm
integral equations of the first kind, when the integral operator
is self--adjoint and compact, is considered here. The data
function is assumed to be perturbed \emph{gently} by an additive
noise so that it still belongs to the range of the operator. First
we estimate upper and lower bounds for the $\varepsilon$--capacity
(and then for the \emph{metric information}), and explicit
computations in some specific cases are given; then the problem is
reformulated from a probabilistic viewpoint and use is made of the
probabilistic information theory. The results obtained by these
two approaches are then compared.
\end{abstract}

\maketitle

\section{Introduction}
\label{se:introduction}
Let us consider the following class of Fredholm integral equations of the first kind:
\beq
\label{1}
Af=g,
\eeq
where $A:X \to Y$ is a self--adjoint compact operator, $X$ and $Y$ being the solution and the data space,
respectively. Hereafter we set $X = Y = L^2[a,b]$.

Solving Equation (\ref{1}) presents two problems:
\begin{itemize}
\item[a)] The $\Ra(A)$ is not closed in the data space $Y$. Therefore, given an arbitrary function
$g\in Y$, it does not follow necessarily that there exists a solution $f \in X$.
\item[b)] Even if two data functions $g_1$ and $g_2$ belong to $\Ra(A)$, and their distance in $Y$ is small,
nevertheless the distance between $A^{-1}g_1$ and $A^{-1}g_2$ can be unlimitedly large, in view of the fact
that the inverse of the compact operator $A$ is not bounded ($X$ and $Y$ being infinite dimensional space).
\end{itemize}
In the numerical applications, $g$ is perturbed by a noise $n$ which can represent either round--off
numerical error or measurement error if $g$ describes experimental data. Assuming in both cases that the
perturbation produced by the noise is additive, the data function actually known is $\barg = g+n$
(instead of the noiseless data function $g$). Then, in order to recover $f$ one is forced to use the
so--called \emph{regularization methods}; the literature on these topics is very extensive, and we shall
return later on this point.

Since the operator $A$ is self--adjoint it admits a set of eigenfunctions $\{\psi_k\}_1^\infty$ and,
accordingly, a countably infinite set of eigenvalues $\{\lambda_k\}_1^\infty$. The eigenfunctions form an
orthonormal basis of the orthogonal complement of the null space of the operator $A$, and therefore an
orthonormal basis of $L^2[a,b]$ when $A$ is injective. For the sake of simplicity we consider hereafter
only this case. The Hilbert--Schmidt theorem guarantees that $\lim_{k\to\infty}\lambda_k = 0$. We shall
suppose hereafter that the eigenvalues are ordered as follows: $\lambda_1 > \lambda_2 > \lambda_3 > \cdots$;
furthermore, we assume for simplicity that they are bounded by 1, i.e., $\lambda_1 \leqslant 1$.
If we consider the noiseless data function $g$, we can associate to the integral equation (\ref{1})
the eigenfunction expansion
\beq
\label{2}
f(x) = \sum_{k=1}^\infty \frac{g_k}{\lambda_k}\,\psi_k(x),
\eeq
where $g_k=(g,\psi_k)$, $(\cdot,\cdot)$ denoting the scalar product in $L^2[a,b]$. The series (\ref{2})
converges in the sense of the $L^2$--norm. Unfortunately this series is not useful since, in practice,
the noiseless data function $g$ is unknown. If we take into account the additive noise $n$, instead of
Equation (\ref{1}), we have
\beq
\label{noisy}
Af + n = \barg.
\eeq
Therefore, instead of expansion (\ref{2}), we have to deal with an expansion of the type
\beq
\label{exp_noisy}
\sum_{k=1}^\infty \frac{\barg_k}{\lambda_k}\,\psi_k(x), \qquad \barg_k=(\barg,\psi_k),
\eeq
which either diverges if $\barg\,\not\in\Ra(A)$, or converges to a function whose distance in norm from the true
solution $f$ (corresponding to the noiseless data) can be quite large. One is then forced to use regularization
procedures as mentioned above.

The mathematical framework outlined so far is only a schematization of reality; in particular,
if the data $g$ describes experimental data, then it obviously will be an element of a finite dimensional
space, while the solution $f$ can still be considered an element of an infinite--dimensional function space;
in general, the data space $Y$ and the solution space $X$ may differ.
In this case the analysis would require the use of singular values and singular functions of the operator $A$
\cite{Bertero,Nashed}, instead of the eigenvalues $\lambda_k$ and eigenvectors $\psi_k$.
For the sake of clarity, here it is convenient to identify data with an element $g$ of $L_2[a,b]$ and deal with
a self--adjoint operator $A$; in this way the analysis is technically simpler, and becomes more transparent
for our purposes.

Several methods of regularization have been proposed \cite{Engl,Groetsch}:
all of them modify one of the elements of the triplet $\{A,X,Y\}$ \cite{Nashed}.
Among these methods, the procedure which is probably the most popular consists in admitting only those solutions
which belong to a compact subset of the solution space $X$. The key theorem used in this method reads as follows:
let $\sigma$ be a continuous map from a compact topological space into a Hausdorff topological space; if
$\sigma$ is one--to--one, then its inverse map $\sigma^{-1}$ is continuous \cite{Kelley}. The condition of
compactness can be realized by the use of \emph{a--priori} bounds \cite{John,Tikhonov}, which require some prior
knowledge or some constraints on the solution. Then the procedure works by taking into account two bounds, one on the
solutions and one on the noise $n$:
\begin{eqnarray}
\label{3a}
\|Bf\|_X & \leqslant & 1,  \\
\label{3b}
\|n\|_Y & \leqslant & \varepsilon,
\end{eqnarray}
where $B$ is a suitable constraint operator. Let us suppose that the eigenfunctions $\{\psi_k\}_1^\infty$
diagonalize the operator $B^*B$, i.e., $A^*A$ and $B^*B$ commute. In such a case we have
$B^*Bf=\sum_{k=1}^\infty\beta_k^2 f_k\psi_k$, where
$f_k=(f,\psi_k)$, and $\beta_k^2$ are the eigenvalues of $B^*B$. The constraint operator $B$ has compact inverse
if and only if $\lim_{k\to\infty}\beta_k^2=+\infty$; under such a condition, the solution obtained by truncating expansion
(\ref{exp_noisy}) at the largest integer $k$ such that $\lambda_k \geqslant \varepsilon\beta_k$, converges to the solution
$f$, as $\varepsilon\to 0$, in the sense of the $L^2$--norm. In several cases a much milder constraint can be
conveniently used, i.e., $B=I$ ($\forall k,\,\beta_k=1$).
In this case the compactness condition, required by the theorem quoted above, is not satisfied; however,
we shall prove in Section \ref{se:metric} that the approximation $f_*$ obtained by truncating expansion
(\ref{exp_noisy}) at the largest $k$ such that $\lambda_k\geqslant\varepsilon$ is convergent, though in weak sense,
to the solution $f$ as $\varepsilon\to 0$.

Hereafter we shall only consider this last truncation method, and we denote by $k_0(\varepsilon)$ the largest
integer $k$ such that $\lambda_k\geqslant\varepsilon$; further, we assume that $\barg\in\Ra(A)$. Since $A$ is compact,
$Y_0 \equiv \Ra(A)$ is a compact subset of $Y$, and then finite coverings of $Y_0$ can be constructed.
By adopting the language of the communication theory \cite{Middleton}, and regarding the inverse problem
of approximating $f$ from a given $\barg$ as a \emph{communication channel} problem, one can compute the
\emph{maximal length} of the messages conveyed back from $\barg$ to $f$.
We are thus led to find a relationship between the maximal length of these messages, which is related to
the truncation number $k_0(\varepsilon)$, and the \emph{massiveness} (or \emph{degree of compactness}) of the set $Y_0$.
It turns out that the \emph{degree of compactness} of $Y_0$ is related to the \emph{smoothness} of the kernel
of the integral operator $A$. In fact, the asymptotic behavior of the eigenvalues $\lambda_k$,
for large $k$, is strictly related to the regularity properties of the kernel: Hille and Tamarkin \cite{Hille}
have systematically explored the relationship between the regularity properties of the kernel and the distribution
of the eigenvalues of the Fredholm integral equation of the first kind. We can say that as the regularity of the
kernel increases, e.g. passing from the class of functions $C^0$ to $C^\infty$ and then to the class of
analytic functions, the eigenvalues $\lambda_k$ decrease more and more rapidly for $k\to\infty$.
Thus the minimum number of balls in a covering of $Y_0$, or the maximum number of balls in a packing of
$Y_0$ \cite{Rogers}, which give a numerical estimate of the \emph{degree of compactness} of $Y_0$, decreases
as the \emph{smoothness} of the kernel increases. Finally, the type of restored continuity in reconstructing
$f$ from a given $\barg$ depends on the \emph{a priori} global bounds imposed on the solution (see formula
(\ref{3a})), and also on the \emph{degree of compactness} of $Y_0$ and, accordingly,
it is related to the length of the messages conveyed back from $\barg$ to reconstruct $f$.
Since we are concerned with the maximal length of these messages we are led to consider a weak--type
convergence in the reconstruction of the solution $f$; accordingly we will define $k_0(\varepsilon)$
as the largest integer such that $\lambda_k \geqslant \varepsilon$. By adopting a more restrictive constraint
we could achieve strong--type convergence, but at the same time we would have shorter messages
conveyed back from $\barg$ for reconstructing $f$.

The problem of reconstructing $f$ from $\barg$ can be reformulated as well in probabilistic terms,
in view of the fact that the data function $g$ is perturbed by the noise $n$, which can be properly regarded
as a random variable. With this in mind one can rewrite equation (\ref{noisy}) in probabilistic form as
\beq
\label{4}
A\xi+\zeta=\eta,
\eeq
where $\xi$, $\zeta$ and $\eta$, which correspond to $f$, $n$ and $\barg$ respectively, are Gaussian weak
random variables \cite{Balakrishnan} in the Hilbert space $L^2[a,b]$.
Next, Equation (\ref{4}) can be turned into an infinite sequence of one--dimensional equations by means of
orthogonal projections, i.e.,
\beq
\label{5}
\lambda_k\xi_k+\zeta_k=\eta_k, \qquad k=1,2,\ldots,
\eeq
where $\xi_k=(\xi,\psi_k)$, $\zeta_k=(\zeta,\psi_k)$, $\eta_k=(\eta,\psi_k)$ are Gaussian random variables.
Using this approach it is possible to evaluate the amount of information $J(\xi_k,\eta_k)$ about the variable $\xi_k$,
which is contained in the variable $\eta_k$. From this approach then another method of truncation emerges,
which is based on neglecting all those components for which $J(\xi_k,\eta_k)$ is less than $\frac{1}{2}\ln 2$.
As illustrated in Section \ref{se:probabilistic}, this criterion leads to a truncation number which is very close
to the number $k_0(\varepsilon)$ introduced previously. One can thus conclude that the two procedures,
the deterministic one, based on the evaluation of the maximal length of the messages conveyed back from $\barg$ to
$f$, and the probabilistic one, based on the information theory, yield essentially the same result.

Information theory, or the theory of coding arose from the fundamental paper of Shannon in 1948 \cite{Shannon}.
Perhaps it should be more correctly referred to as statistical communication theory. The information source is any
producer of information according to some known probability law, and this information has to be
communicated to the destination by means of a transmission channel. Noise can be regarded as anything which
impairs the ability of the channel to transmit with complete reliability. Information theory is concerned with the
methods for achieving high reliability without reducing the transmission rate too drastically. Successively
the mathematical theory of information was extended by several authors, notably Kolmogorov and Gelfand
(see, in particular, \cite{Gelfand1} and the papers quoted therein). One question quite naturally arises:
On the one hand information theory is formulated in the framework and uses language and tools of the probability
theory, on the other hand the concept of \emph{information} can be thought of as more basic and independent
of probability \cite{Kolmogorov2}.
Then the problem becomes: how to construct a nonprobabilistic theory of information.
To this purpose Kolmogorov and his school introduced and developed an alternative approach to the quantitative
definition of information, which is logically independent of probabilistic assumptions: the measure of
information is given in purely combinatorial terms \cite{Kolmogorov2}.
This combinatorial, or metric, approach finally results in the theory of the $\varepsilon$--entropy and
$\varepsilon$--capacity of sets in metric spaces \cite{Kolmogorov1}.

The connection between ideas and concepts of Shannon's information theory, with particular attention to the
notion of \emph{length of a message in binary units}, and those of $\varepsilon$--entropy and $\varepsilon$--capacity
are illustrated in detail in \cite{Kolmogorov1}, to which the interested reader is referred
(to this purpose, let us also mention \cite{VandeGeer}, where the $\varepsilon$--entropy plays a
crucial role in connection with empirical processes estimation).
With a small abuse of language we call \emph{metric information} that induced by the $\varepsilon$--capacity,
which is, indeed, defined as the number of binary signs that can be reliably transmitted.
Finally, the problem of comparing the results of probabilistic and nonprobabilistic, or metric, information theory
remains. The main aim of this paper consists precisely in trying to give a partial answer to this question
in the specific case of Fredholm integral equations of the first kind.

The paper is organized as follows. In Section \ref{se:metric} we first prove that the approximation $f_*$
converges weakly to $f$ as $\varepsilon\to 0$. Then we find an upper and a lower bound for the $\varepsilon$--entropy
associated with the mapping of the unit ball, in the solution space, induced by the operator $A$.
Next, we evaluate explicitly an upper bound for the maximal length of the messages conveyed
back from $\barg$ to reconstruct $f$, and this provides an estimate of what we call \emph{metric information}.
Explicit calculations are given in three specific cases: harmonic continuation, backward solution of the
heat equation, first kind Fredholm integral equation with continuous kernel. In Section \ref{se:probabilistic}
we reconsider the problem from a probabilistic viewpoint. We introduce another truncation method based on
probabilistic information theory, and accordingly we derive an approximation which converges to the solution,
in the sense of the probabilistic theory, under suitable conditions on the covariance operator of the solution.

\section{Metric Information Associated with Fredholm Integral Equations of the First Kind}
\label{se:metric}
\subsection{Weak convergence of the $\mathbf{f_*}$ approximation}
\label{subse:weak}
Let us consider the approximation $f_*=\sum_{k=1}^{k_0(\varepsilon)}(\barg_k/\lambda_k)\psi_k$ where $k_0(\varepsilon)$ is
the largest integer such that $\lambda_k\geqslant\varepsilon$.
We want to prove the weak convergence of $f_*$ to $f$ as $\varepsilon\to 0$ and, accordingly, the \emph{weak continuity}
in the restored solution; for this purpose we need the following auxiliary lemma.

\begin{lemma}
\label{weaklemma}
For any function $f$ which satisfies the following bounds
\begin{eqnarray}
\label{w1}
\left\|Af-\barg\right\|_{Y \equiv L^2[a,b]} &\leqslant & \varepsilon, \\
\label{w2}
\left\|f\right\|_{X \equiv L^2[a,b]} &\leqslant & 1,
\end{eqnarray}
the following inequalities hold:
\begin{eqnarray}
\left\|A(f-f_*)\right\|_Y &\leqslant& \sqrt{2}\varepsilon,\label{w3} \\
\left\|f-f_*\right\|_X &\leqslant& \sqrt{2},\label{w4} \\
\left\|A(f-f_*)\right\|_Y^2 + \varepsilon^2 \left\|f-f_*\right\|_X^2 &\leqslant& 4\varepsilon^2.\label{w5}
\end{eqnarray}
\end{lemma}

\begin{proof}
(a) From the inequality $\lambda_k < \varepsilon$ for $k>k_0$ and the bound $\|f\|_X \leqslant 1$ it follows:
\beq
\label{w6}
\sum_{k=k_0+1}^\infty \lambda_k^2 |f_k|^2 < \varepsilon^2.
\eeq
From $\|Af-\barg\|_Y\leqslant\varepsilon$ we get:
\beq
\label{w7}
\sum_{k=1}^{k_0} \lambda_k^2 \left|f_k-\frac{\barg_k}{\lambda_k}\right|^2 \leqslant \varepsilon^2.
\eeq
Therefore we have
\beq
\label{w8}
\left\|A(f-f_*)\right\|^2_Y = \sum_{k=1}^{k_0} \lambda_k^2 \left|f_k-\frac{\barg_k}{\lambda_k}\right|^2
+ \sum_{k=k_0+1}^\infty \lambda_k^2 |f_k|^2 \leqslant 2\varepsilon^2,
\eeq
and inequality (\ref{w3}) is proved. \\
(b) From the inequality $\lambda_k\geqslant\varepsilon$ for $k\leqslant k_0$ and the bound $\|Af-\barg\|_Y\leqslant\varepsilon$ we obtain
\beq
\label{w9}
\sum_{k=1}^{k_0} \left|f_k-\frac{\barg_k}{\lambda_k}\right|^2 =
\sum_{k=1}^{k_0} \frac{1}{\lambda_k^2} \left|\lambda_k f_k-\barg_k\right|^2 \leqslant 1.
\eeq
From $\|f\|_X\leqslant 1$ it follows:
\beq
\label{w10}
\sum_{k=k_0+1}^\infty |f_k|^2 \leqslant 1.
\eeq
Therefore we have:
\beq
\label{w11}
\left\|f-f_*\right\|_X^2= \sum_{k=1}^{k_0} \left|f_k-\frac{\barg_k}{\lambda_k}\right|^2 +
\sum_{k=k_0+1}^\infty |f_k|^2 \leqslant 2,
\eeq
and inequality (\ref{w4}) is proved. Next, from (\ref{w8}) and (\ref{w11}) we obtain:
\beq
\label{w12}
\left\|A(f-f_*)\right\|^2_Y + \varepsilon^2 \left\|f-f_*\right\|^2_X \leqslant 4\varepsilon^2,
\eeq
that is, inequality (\ref{w5}).
\end{proof}

Let us note that $\lim_{\varepsilon\to 0}k_0(\varepsilon)=+\infty$. The latter equality follows from the
definition itself of $k_0(\varepsilon)$ and from the fact that $\lim_{k\to\infty}\lambda_k=0$. Next
we prove the following theorem.

\begin{theorem}
\label{theo1}
For any function $f$ which satisfies bounds $(\ref{w1})$ and $(\ref{w2})$, the following limit holds true:
\beq
\label{w13}
\lim_{\varepsilon\to 0} \left(f-f_*,v\right)_X = 0, \qquad \forall v\in X;\,\|v\|_X\leqslant 1.
\eeq
\end{theorem}

\begin{proof}
Let us put: $x_k = f_k-(f_*)_k$; then we have:
\beq
\label{w14}
\left(f-f_*,v\right)_X = \sum_{k=1}^\infty x_k v_k, \qquad \left(\sum_{k=1}^\infty |v_k|^2 \leqslant 1\right).
\eeq
Next, by the Schwarz inequality and bound (\ref{w5}), we have:
\beq
\label{w15}
\begin{split}
\sds\left|\left(f-f_*,v\right)_X\right| &\leqslant \sds\sum_{k=1}^\infty |x_k v_k| =
\sum_{k=1}^\infty \left(\frac{\lambda_k^2+\varepsilon^2}{\lambda_k^2+\varepsilon^2} \right)^{1/2}|x_k v_k| \\
&\leqslant \left(\left\|A(f-f_*)\right\|^2_Y + \varepsilon^2 \left\|f-f_*\right\|^2_X
\right)^{1/2}\left(\sum_{k=1}^\infty\frac{|v_k|^2}{\lambda_k^2+\varepsilon^2}\right)^{1/2} \\
&\leqslant\left(4\varepsilon^2\sum_{k=1}^\infty\frac{|v_k|^2}{\lambda_k^2+\varepsilon^2}\right)^{1/2}.
\end{split}
\eeq
Next we split the sum $\sum_{k=1}^\infty |v_k|^2/(\lambda_k^2+\varepsilon^2)$ into two parts, i.e.,
\beq
\label{w16}
\sum_{k=1}^{k_0}\frac{|v_k|^2}{\lambda_k^2+\varepsilon^2} +
\sum_{k=k_0+1}^\infty \frac{|v_k|^2}{\lambda_k^2+\varepsilon^2}.
\eeq
The first term of the sum (\ref{w16}) can be majorized as follows:
\beq
\label{w17}
\sum_{k=1}^{k_0}\frac{|v_k|^2}{\lambda_k^2+\varepsilon^2} \leqslant \frac{1}{2\varepsilon^2} \sum_{k=1}^\infty |v_k|^2
\leqslant \frac{1}{2\varepsilon^2}.
\eeq
From formulae (\ref{w15}) and (\ref{w17}) we have
\beq
\label{w18}
4\varepsilon^2\sum_{k=1}^{k_0}\frac{|v_k|^2}{\lambda_k^2+\varepsilon^2} \leqslant
2 \sum_{k=1}^\infty |v_k|^2 \leqslant 2.
\eeq
Furthermore, $\lim_{\varepsilon\to 0}\varepsilon^2|v_k|^2/(\lambda_k^2+\varepsilon^2)=0$ for $k\leqslant k_0$.
Therefore we have
\beq
\label{w19}
\lim_{\varepsilon\to 0} 4\varepsilon^2\sum_{k=1}^{k_0}\frac{|v_k|^2}{\lambda_k^2+\varepsilon^2}=0.
\eeq
Let us now consider the second term of sum (\ref{w16}); we can write
\beq
\label{w20}
\sum_{k=k_0+1}^\infty \frac{|v_k|^2}{\lambda_k^2+\varepsilon^2}\leqslant
\frac{1}{\varepsilon^2}\sum_{k=k_0+1}^\infty |v_k|^2.
\eeq
Therefore from formulae (\ref{w15}) and (\ref{w20}) we get
\beq
\label{w21}
4\varepsilon^2\sum_{k=k_0+1}^\infty \frac{|v_k|^2}{\lambda_k^2+\varepsilon^2}\leqslant
4\sum_{k=k_0+1}^\infty |v_k|^2.
\eeq
Then, taking into account that $\lim_{\varepsilon\to 0} k_0(\varepsilon)=+\infty$, we can conclude:
\beq
\label{w22}
\lim_{\varepsilon\to 0} 4\varepsilon^2\sum_{k=k_0+1}^\infty \frac{|v_k|^2}{\lambda_k^2+\varepsilon^2}=0.
\eeq
From (\ref{w19}) and (\ref{w22}) we then obtain:
\beq
\label{w23}
\lim_{\varepsilon\to 0} \left(f-f_*,v\right)_X=0, \qquad \forall v\in X;\,\|v\|_X\leqslant 1,
\eeq
and the theorem is proved.
\end{proof}

\subsection{$\mathbf{\varepsilon}$--entropy and $\mathbf{\varepsilon}$--capacity associated with the operator $\mathbf{A}$}
\label{subse:epsilon}
Let us consider the unit ball in the solution space $X \equiv L^2[a,b]$, i.e., $\{f \in X \,| \, \|f\|_{X} \leqslant 1\}$.
The operator $A$ maps the unit ball onto a compact ellipsoid $\cE \in \Ra(A)$ contained in $Y \equiv L^2[a,b]$,
whose semi--axes' lengths are the eigenvalues $\lambda_k$ of the operator $A$.
In order to give a numerical estimate of the \emph{massiveness} of the set $\cE$, let us first recall some
basic definitions \cite{Kelley,Lorentz}:
\begin{itemize}
\item[(a)] A family $Y_1,\cdots,Y_n$ of subsets of $Y$ is an $\varepsilon$--covering of $\cE$ if the diameter
of each $Y_k$ does not exceed $2\varepsilon$ and if the sets $Y_k$ cover $\cE$: $\cE \subset\cup_{k=1}^n Y_k$.
\item[(b)] Points $y_1,\cdots,y_m$ of $\cE$ are called $\varepsilon$--distinguishable if the distance between
each two of them exceeds $\varepsilon$.
\end{itemize}
Since $\cE$ is compact, then there exists a finite $\varepsilon$--covering for each $\varepsilon >0$ and, moreover,
$\cE$ can contain only finitely many $\varepsilon$--distinguishable points.
For a given $\varepsilon>0$, the number $n$ of sets $Y_k$ in a covering family depends on the family, but the
minimal value of $n$, $N_\varepsilon(\cE)=\min n$, is an invariant of the set $\cE$, which depends only on $\varepsilon$.
Its logarithm (throughout the paper $\log x$ will always denote the logarithm of the number $x$ to
the base 2), that is, the function $H_\varepsilon (\cE)=\log N_\varepsilon(\cE)$ is the $\varepsilon$--entropy
of the set $\cE$.
Analogously, the number $m$ in definition (b) depends on the choice of points, but its maximum
$M_\varepsilon(\cE)=\max m$ is an invariant of the set $\cE$. Its logarithm, that is the function
$C_\varepsilon(\cE)=\log M_\varepsilon(\cE)$ is called the $\varepsilon$--capacity of the set $\cE$. This quantity
represents the maximum number of $\varepsilon$--distinguishable signals that can be received, that is those
data which satisfy the following inequalities $\|\barg^{(i)}-\barg^{(k)}\|_{Y} > \varepsilon$,
for all $i \neq k$, $\barg^{(i)}$, $\barg^{(k)} \in \cE$.

A general result about $\varepsilon$--entropy and $\varepsilon$--capacity are the following inequalities
\cite{Lorentz}:
\beq
\label{7}
H_\varepsilon (\cE) \leqslant C_\varepsilon(\cE) \leqslant H_{\varepsilon/2} (\cE).
\eeq
To obtain estimates for the $\varepsilon$--capacity $C_\varepsilon(\cE)$, our aim
now is to look for a lower bound for $H_\varepsilon (\cE)$ and an
upper bound for $H_{\varepsilon/2} (\cE)$. For this purpose, let us
consider the finite dimensional subspace $Y_{k_0}$ of $Y$ spanned
by the first $k_0$ axes of $\cE$, and put $\cE_{k_0} = \cE \cap
Y_{k_0}$. Then $\cE_{k_0}$ is a finite dimensional ellipsoid whose
volume is just $\prod_{k=1}^{k_0}\lambda_k$ times the volume
$\Omega_{k_0}$ of the unit ball in $Y_{k_0}$. Since the volume of
an $\varepsilon$--ball in $Y_{k_0}$ is just $\varepsilon^{k_0}
\Omega_{k_0}$, we see that in order to cover the ellipsoid $\cE$
by $\varepsilon$--balls we shall need at least
$\prod_{k=1}^{k_0}\lambda_k/\varepsilon$ such balls. From this it
follows that:
\beq
\label{8}
\prod_{k=1}^{k_0}\frac{\lambda_k}{\varepsilon} \leqslant N_\varepsilon(\cE),
\eeq
and therefore we have the following lower bound for the
$\varepsilon$--entropy $H_\varepsilon (\cE)$:
\beq \label{9}
\sum_{k=1}^{k_0} \log \frac{\lambda_k}{\varepsilon} \leqslant \log N_\varepsilon(\cE) = H_\varepsilon (\cE).
\eeq
An upper bound for $H_{\varepsilon/2}(\cE)$ can be found in the following way
\cite{Gelfand2,Prosser}: Let us construct in $Y_{k_0}$ the cubical
lattice with mesh width $\varepsilon_1=\varepsilon/(2\sqrt{k_0})$, and
with coordinate axes the axes of $\cE_{k_0}$. In view of the
choice of $\varepsilon_1$ any point of $Y_{k_0}$, and in particular
of $\cE_{k_0}$, lies within a distance not exceeding
$\frac{1}{2}\varepsilon_1\sqrt{k_0}=(\varepsilon/4)$ from the
nearest point of this lattice. In particular, it will lie at a
distance not exceeding $(\varepsilon/4)$ from one of the
lattice points which are contained in the parallelepiped $P_{k_0}$
defined by:
\beq
\label{parallelpiped}
-\frac{\varepsilon}{4}-\lambda_k \leqslant x_k \leqslant
\frac{\varepsilon}{4}+\lambda_k, \qquad 1\leqslant k \leqslant k_0.
\eeq
Now, if $k_0=k_0(\varepsilon/4)$, that is $k_0$ represents the number
of terms in the sequence $\{\lambda_k\}$ which are greater than
$(\varepsilon/4)$, then every point $x\in\cE$ lies at a
distance not exceeding $(\varepsilon/4)$ from a point of
$\cE_{k_0}$. In fact, let us write $x=\sum_k x_k \psi_k$,
$\{\psi_k\}$ being the orthonormal basis for $Y$ made of the
eigenvectors of the operator $A$. Since $x$ belongs to $\cE$, then
evidently
$\sum_{k=1}^\infty\left|x_k/\lambda_k\right|^2\leqslant 1$.
Hence the square of the distance from $x$ to $\cE_{k_0}$ is
\beq
\label{distance}
d^2(x,\cE_{k_0})=\sum_{k=k_0+1}^\infty|x_k|^2=\sum_{k=k_0+1}^\infty\lambda_k^2\left|\frac{x_k}{\lambda_k}\right|^2
\leqslant \lambda_{k_0+1}^2
\sum_{k=1}^\infty\left|\frac{x_k}{\lambda_k}\right|^2 \leqslant
\left(\frac{\varepsilon}{4}\right)^2.
\eeq
Now, the balls of radius $(\varepsilon/2)$ with centers at those lattice points within
$P_{k_0}$ cover the ellipsoid $\cE$. In fact, from
(\ref{distance}) each point of $\cE$ is at a distance not
exceeding $(\varepsilon/4)$ from $\cE_{k_0}$, and each point of
$\cE_{k_0}$ is at a distance not exceeding $(\varepsilon/4)$
from some point of the lattice belonging to $P_{k_0}$; then each
point of $\cE$ lies at a distance not exceeding
$(\varepsilon/2)$ from some point of the lattice belonging to
$P_{k_0}$. Obviously the number of lattice points in $P_{k_0}$ is
not greater than
\beq
\label{number}
\prod_{k=1}^{k_0} 2\left(\frac{\lambda_k}{\varepsilon_1}+1\right)=
\prod_{k=1}^{k_0}\frac{2}{\varepsilon}\,(2\lambda_k\sqrt{k_0}+\varepsilon)\leqslant
\left(\frac{6\,\sqrt{k_0}}{\varepsilon}\right)^{k_0},
\eeq
where we used the assumption $\varepsilon<\lambda_1\leqslant 1\leqslant k_0$. Then the
number of elements in this $\varepsilon$--covering is no more than
$\left[6\,\sqrt{k_0(\varepsilon/4)}/\varepsilon\right]^{k_0(\varepsilon/4)}$
since $k_0=k_0(\varepsilon/4)$. Taking the logarithm, we finally
obtain
\beq
\label{10}
H_{\varepsilon/2} (\cE) \leqslant
k_0\left(\frac{\varepsilon}{4}\right)\log\frac{6\sqrt{k_0(\varepsilon/4)}}{\varepsilon}
=
k_0\left(\frac{\varepsilon}{4}\right)\left[\log\left(\frac{1}{\varepsilon}\right)+\log
6 + \frac{1}{2}\log k_0\left(\frac{\varepsilon}{4}\right)\right].
\eeq
For the next step we note that $H_\varepsilon(\cE)$ is a
nondecreasing function as $\varepsilon\to 0$, then we can introduce
the \emph{order of growth} $\rho(\cE)$ of the entropy $H_\varepsilon
(\cE)$ as follows:
\beq
\label{11}
\rho(\cE) = \lim_{\varepsilon\to 0}\sup\frac{\log H_\varepsilon(\cE)}{\log(1/\varepsilon)},
\eeq
or, in the case $\rho(\cE)=0$, the \emph{logarithmic order of growth}
$\sigma(\cE)$ of $H_\varepsilon(\cE)$ which reads
\beq
\label{12}
\sigma(\cE)= \lim_{\varepsilon\to 0}\sup\frac{\log H_\varepsilon(\cE)}{\log\log(1/\varepsilon)}.
\eeq
Since we are interested in relating the asymptotic behavior of
$H_\varepsilon(\cE)$ as $\varepsilon\to 0$ with the asymptotic behavior
of the \emph{semi--axes} $\{\lambda_k\}$ of $\cE$ as $k\to\infty$,
we are led to introduce the \emph{exponent of convergence}
$\lambda$ and the \emph{logarithmic exponent of convergence} $\mu$
of the sequence $\{1/\lambda_k\}$, see \cite{Levin}:
\begin{eqnarray}
\lambda & = & \lim_{\varepsilon\to 0}\sup\frac{\log k_0(\varepsilon)}{\log(1/\varepsilon)}, \label{13}\\
\mu & = & \lim_{\varepsilon\to 0}\sup\frac{\log k_0(\varepsilon)}{\log\log(1/\varepsilon)},\label{14}
\end{eqnarray}
where $k_0(\varepsilon)$ denotes the number of elements of the sequence ${\lambda_k}$
which are greater than $\varepsilon$.
The following relationship is proved in \cite{Prosser}: $\rho(\cE)=\lambda$, and if
$\rho(\cE)=\lambda=0$, then $\sigma(\cE)=\mu+1$. Finally, we can define the \emph{degree of compactness}
$d_c$ associated with the range of the operator $A$ as $d_c = 1/\rho$ (if $\rho \neq 0$),
and the \emph{exponential degree of compactness} of $\Ra(A)$ as $d_c^{\,e} = 2^{1/\sigma}$ (if $\rho=0$).

By using bounds (\ref{9}) and (\ref{10}), we can now evaluate the degree of compactness of $\Ra(A)$
in three specific examples: harmonic continuation, backward solution of the heat equation, and a convolution
equation with continuous kernel; in all these examples the behavior with $k$ of the eigenvalues is uniform,
in the sense that the relative rate of decaying of the eigenvalues follows, for all $k$, a uniform law in $k$.

\subsubsection{Harmonic continuation}
\label{subsubse:harmonic}
Let us consider a family $\cF$ of functions $u(r,\theta)$ which satisfy the Laplace equation at the interior
of the unit disk. We want to determine $u(b,\theta)$, ($b<1$), assuming that $u(a,\theta)$ ($a<b$) is known
within a certain approximation. The solution to the problem is obtained by solving the following integral
equation of Fredholm--type:
\beq
\label{15}
u(a,\theta)=\frac{1}{2\pi}\int_0^\pi P(\theta-\phi)\,u(b,\phi)\,d\phi, \qquad -\pi < \theta\leqslant\pi,
\eeq
where $P(\theta-\phi)$ is the Poisson kernel given by:
\beq
\label{16}
P(\theta-\phi) = \sum_{k=-\infty}^{+\infty}\left(\frac{a}{b}\right)^{|k|}e^{ik(\theta-\phi)}.
\eeq
We can put Equation (\ref{15}) into the form (\ref{1}): $Af=g$, where $f(\phi) \equiv u(b,\phi)$,
$g(\theta)=u(a,\theta)$, $(b>a)$; $u(b,\phi)$ is the restriction to the circle of radius $b$ of a function
harmonic in the unit disk, which belongs to $L^2[-\pi,\pi]$; then the following expansion converges in the sense of the
$L^2$--norm:
\beq
\label{17}
u(1,\theta)=\sum_{k=-\infty}^{+\infty} u_k e^{ik\theta}, \qquad \left(\sum_{k=-\infty}^{+\infty}|u_k|^2 < \infty\right).
\eeq
Furthermore, we have:
\beq
\label{18}
u(b,\theta)=\sum_{k=-\infty}^{+\infty} b^{|k|} u_k e^{ik\theta},
\eeq
which is uniformly convergent.
The eigenvalues of the operator $A$ are $\lambda_k=(a/b)^{|k|}$, $b>a$, and the eigenfunctions are given
by $\psi_k(\theta) = e^{-ik\theta}$; evidently, $\lim_{k\to\infty}\lambda_k=0$.
The $\Ra(A)$ is not closed in $L^2[-\pi,\pi]$; in fact, only those functions $u$ which satisfy the following bound:
\beq
\label{19}
\sum_{k=-\infty}^{+\infty} \left(u_k a^{|k|}\right)^2 < \infty,
\eeq
belong to the $\Ra(A)$. Now, if a noise $n$ is added to the data function $g$, the function actually
known is $\barg=g+n$ which, in general, does not belong to $\Ra(A)$; nevertheless hereafter we still assume that
$\barg \in \Ra(A)$. Next we restrict the solution space to those functions which satisfy the following bound:
\beq
\label{20}
\sum_{k=-\infty}^{+\infty} \left(u_k b^{|k|}\right)^2 \leqslant 1.
\eeq
It is now easy to evaluate the truncation number $k_0(\varepsilon)$, which is given by the largest integer such
that $\lambda_k \geqslant \varepsilon$, i.e.,
\beq
\label{21}
k_0(\varepsilon) = \left[\frac{\log(1/\varepsilon)}{\log(b/a)}\right],
\eeq
where $[\cdot]$ stands for the integral part.
Now we split the sums (\ref{16})-(\ref{20}) into two parts: the first is obtained by varying $k$ from zero to $+\infty$;
the second by varying $k$ from $-1$ to $-\infty$. We denote the $\varepsilon$--entropy ($\varepsilon$--capacity) associated
with the truncation of the first sum by $H_\varepsilon^{(+)}(\cE)$ ($C_\varepsilon^{(+)}(\cE)$);
accordingly, the $\varepsilon$--entropy ($\varepsilon$--capacity) associated with the truncation of the second sum
by $H_\varepsilon^{(-)}(\cE)$ ($C_\varepsilon^{(-)}(\cE)$).
Then using formula (\ref{21}) and inequality (\ref{10}) we obtain:
\beq
\label{22}
\begin{split}
\sum_{k=1}^{k_0(\varepsilon)}\log\left(\frac{\lambda_k}{\varepsilon}\right) &\leqslant
H_\varepsilon^{(+)}(\cE) \leqslant C_\varepsilon^{(+)}(\cE) \\
&\leqslant H_{\varepsilon/2}^{(+)}(\cE)\leqslant
k_0\left(\frac{\varepsilon}{4}\right)
\left[\log\left(\frac{1}{\varepsilon}\right)+\log 6+
\frac{1}{2}\log k_0\left(\frac{\varepsilon}{4}\right)\right] \\
&\leqslant \frac{2+\log(1/\varepsilon)} {\log(b/a)}
\left[\log\left(\frac{1}{\varepsilon}\right)+\log 6+\frac{1}{2}\log k_0\left(\frac{\varepsilon}{4}\right)\right].
\end{split}
\eeq
The leading term on the r.h.s. of (\ref{22}) as $\varepsilon\to 0$ is given by
\beq
\label{leading_1}
\frac{\log(1/\varepsilon)} {\log(b/a)} \log\left(\frac{1}{\varepsilon}\right) \sim
k_0(\varepsilon) \log\left(\frac{1}{\varepsilon}\right),
\eeq
while the leading term on the l.h.s. of (\ref{22}) becomes
\beq
\label{leading_11}
\frac{1}{2} k_0(\varepsilon) \log\left(\frac{1}{\varepsilon}\right).
\eeq
We thus obtain, for $\varepsilon$ sufficiently small, fairly sharp inequalities for the $\varepsilon$--capacity:
\beq
\label{23}
\frac{1}{2} k_0(\varepsilon) \log\left(\frac{1}{\varepsilon}\right)
\lesssim C_\varepsilon^{(+)}(\cE)\lesssim
k_0(\varepsilon)\log\left(\frac{1}{\varepsilon}\right) \leqslant \frac{\log^2(1/\varepsilon)}{\log(b/a)}.
\eeq
We thus have an upper bound for the maximal length, in binary units, of the messages conveyed back from $\barg$
to reconstruct $f$, associated with the truncation of the positive sum; we obtain, with obvious notation:
\beq
\label{24}
L_{\max}^{(+)}(\varepsilon)\lesssim 2^{k_0(\varepsilon)\log(1/\varepsilon)} \sim 2^{(\log^2(1/\varepsilon)/\log(b/a))}.
\eeq
Finally, for the total maximal length we obtain:
\beq
\label{24bis}
\begin{split}
L_{\max}(\varepsilon) &= L_{\max}^{(+)}(\varepsilon)+L_{\max}^{(-)}(\varepsilon)\lesssim
2^{k_0(\varepsilon)\log(1/\varepsilon)+1} \\
&\staccrel{\sim}{\varepsilon\to 0} 2^{k_0(\varepsilon)\log(1/\varepsilon)} \sim
2^{(\log^2(1/\varepsilon)/\log(b/a))},
\end{split}
\eeq
which can be taken as a quantitative estimate of the \emph{metric information}.

\begin{remark}
Let us note that $\log(b/a) = \text{\rm Cons.}\cdot L\{C\}$, where $L\{C\}$ is the extremal length of $\{C\}$, the latter
expressing the set of curves in the ring domain $0<a<r<b<\infty$, which join $r=a$ to $r=b$. $L\{C\}$ is
a conformal invariant \cite{Fuchs}. The r.h.s. of (\ref{23}) may be regarded as a particular case of a more
general result due to Erohin (see \cite{Lorentz}), which shows that for general sets of analytic functions:
\beq
\label{ero}
H_\varepsilon^{(+)} \sim C_\varepsilon^{(+)} \sim \gamma \log^2\left(\frac{1}{\varepsilon}\right),
\eeq
$\gamma$ depending on some conformal invariant.
\end{remark}

Concerning the \emph{order of growth} $\rho(\cE)$ of the $\varepsilon$--entropy and the \emph{exponent of convergence}
$\lambda$: from (\ref{21}) it follows that $\rho(\cE)=\lambda=0$. We then move on to the
\emph{logarithmic order of growth} $\sigma(\cE)$ and, correspondingly, to the
\emph{logarithmic exponent of convergence} $\mu$; we have $\sigma(\cE)=2$ and, consequently,
the \emph{exponential degree of compactness} $d_c^{\,e}=2^{1/\sigma}=2^{1/2}$.

\subsubsection{Backward solution of the heat equation}
\label{subsubse:backward}
Let us consider a heat conducting ring of radius 1. One can pose two problems:
\begin{itemize}
\item[i)] \emph{Direct problem.} Determine the temperature distribution $h(t,\theta)$ at time $t$, when
$h(0,\theta)$ is given. The solution is obtained by solving the Cauchy problem for the heat equation:
\begin{eqnarray}
\label{25}
h_t & = & D\,h_{\theta\theta}, \qquad D > 0, \\
h(0,\theta) & = & h_0(\theta), \qquad 0\leqslant \theta < 2\pi.
\end{eqnarray}
\item[ii)] \emph{Inverse problem.} Determine the temperature distribution $h(b,\theta) = f(\theta)$, at time
$t=b$, when $h(a,\theta) \equiv g(\theta)$, $a>b$, is given.
\end{itemize}
The solution to the inverse problem is obtained by solving the Fredholm integral equation of the first kind:
\beq
\label{27}
h(a,\theta)\equiv g(\theta)=\frac{1}{2\pi}\int_{-\pi}^{\pi}\cK (\theta-\phi) f(\phi)\,d\phi,
\eeq
where the kernel $\cK(\theta-\phi)$ is the elliptic Jacobi theta function:
\beq
\label{28}
\cK (\theta-\phi)=\sum_{k=-\infty}^{+\infty} e^{-Dk^2(a-b)}e^{ik(\theta-\phi)}.
\eeq
The eigenfunctions and the eigenvalues of the integral operator $A$ are respectively $\psi_k(\theta)=e^{-ik\theta}$,
$\lambda_k=\exp(-Dk^2(a-b))$; moreover, $\lim_{k\to\infty}\lambda_k=0$.
Once again we assume that the solution and the data space $X$ and $Y$
are both $L^2[-\pi,\pi]$. We may now consider the following expansion
\beq
\label{29}
h(t,\theta)=\sum_{k=-\infty}^{+\infty} h_k e^{-Dk^2t}e^{ik\theta},
\eeq
which converges in the sense of the $L^2$--norm.

Again the $\Ra(A)$ is not closed in $L^2[-\pi,\pi]$; in fact only those functions $h$ which satisfy the following bound:
\beq
\label{30}
\sum_{k=-\infty}^{+\infty}\left(h_k\,e^{-Dk^2a}\right)^2 < \infty,
\eeq
belong to $\Ra(A)$. If a noise $n$ is added to the data function $g$, only the function
$\barg=g+n$ is known and, in general, it does not belong to $\Ra(A)$. Nevertheless we assume even in this case that
$\barg\in\Ra(A)$. Next we restrict the solution space to a subspace composed of those functions which satisfy the
following a--priori constraint:
\beq
\label{31}
\|h\|^2_{L^2}=\sum_{k=-\infty}^{+\infty}\left(h_k\,e^{-Dk^2b}\right)^2  \leqslant 1.
\eeq
The truncation number $k_0(\varepsilon)$, which is given by the largest integer such that $\lambda_k\geqslant\varepsilon$
can be easily evaluated, i.e.,
\beq
\label{32}
k_0(\varepsilon) = \left[\left(\frac{\log(1/\varepsilon)}{D(a-b)}\right)^{1/2}\right].
\eeq
Based on considerations analogous to those developed in the case of harmonic continuation, and by splitting the sums
(\ref{28})--(\ref{31}) into two sums as done before, we obtain:
\beq
\label{31bis}
\begin{split}
\sum_{k=1}^{k_0(\varepsilon)}\log\left(\frac{\lambda_k}{\varepsilon}\right)&\leqslant C_\varepsilon^{(+)}(\cE) \leqslant
k_0\left(\frac{\varepsilon}{4}\right)\left[\log\left(\frac{1}{\varepsilon}\right)+\log 6+
\frac{1}{2}\log k_0\left(\frac{\varepsilon}{4}\right)\right] \\
&\leqslant \left(\frac{2+\log(1/\varepsilon)} {D(b-a)}\right)^{1/2}
\left[\log\left(\frac{1}{\varepsilon}\right)+\log 6+\frac{1}{2}\log k_0\left(\frac{\varepsilon}{4}\right)\right].
\end{split}
\eeq
The leading term on the r.h.s. of (\ref{31bis}), as $\varepsilon\to 0$, is given by
\beq
\label{leading_2}
\left(\frac{\log(1/\varepsilon)}{D(a-b)}\right)^{1/2} \log\left(\frac{1}{\varepsilon}\right)
\sim k_0(\varepsilon)\log(1/\varepsilon),
\eeq
while the leading term on the l.h.s. of (\ref{31bis}), as $\varepsilon\to 0$, is
\beq
\label{leading_22}
\left(1-\frac{1}{3}\log e\right) k_0(\varepsilon)\log(1/\varepsilon).
\eeq
We therefore have quite sharp bounds on the $\varepsilon$--capacity, i.e.,
\beq
\label{leading_222}
\left(1-\frac{1}{3}\log e\right) k_0(\varepsilon)\log(1/\varepsilon) \lesssim
C_\varepsilon^{(+)}(\cE) \lesssim
k_0(\varepsilon)\log(1/\varepsilon).
\eeq
Then, we have an upper bound for the maximal length, in binary units, of the messages conveyed back
from the data for reconstructing the solution, i.e.,
\beq
\label{33}
L_{\max}^{(+)}(\varepsilon) \lesssim 2^{k_0(\varepsilon)\log(1/\varepsilon)} \leqslant
2^{\,\frac{\text{\rm Cons.}}{(a-b)^{1/2}}[\log(1/\varepsilon)]^{3/2}}.
\eeq
Then the final result referring to the total maximal length of the messages is:
\beq
\label{33bis}
\begin{split}
L_{\max}(\varepsilon) &= L_{\max}^{(+)}(\varepsilon)+L_{\max}^{(-)}(\varepsilon) \lesssim
2^{k_0(\varepsilon)\log(1/\varepsilon)+1} \staccrel{\sim}{\varepsilon\to 0} 2^{k_0(\varepsilon)\log(1/\varepsilon)} \\
&\sim 2^{\,\frac{\smCs}{(a-b)^{1/2}}[\log(1/\varepsilon)]^{3/2}}.
\end{split}
\eeq
Accordingly, the \emph{exponential degree of compactness} is given by $d_c^{\,e}=2^{2/3}$.

\subsubsection{First kind Fredholm integral equation with continuous kernels}
\label{subsubse:first}
Let us consider the following Fredholm integral equation of the first kind:
\beq
\label{34}
Af \equiv \int_0^1 \cK(x,y)\,f(y)\,dy\,= g(x),
\eeq
where the kernel $\cK(x,y)$ is the continuous function
\begin{eqnarray}
\label{35}
\cK(x,y)&=&(1-x)y, \qquad 0\leqslant y\leqslant x \leqslant 1, \\
\label{36}
\cK(x,y)&=&x(1-y), \qquad 0\leqslant x\leqslant y \leqslant 1.
\end{eqnarray}
Eigenfunctions and eigenvalues of operator $A$ in Equation (\ref{34})
can be easily evaluated: the eigenvalues are:
$\lambda_k=1/(k^2\pi^2)$. Once again, following considerations
analogous to those developed in the previous examples we obtain
$k_0(\varepsilon) = [1/(\pi\sqrt{\varepsilon})]$ and, for $\varepsilon$
sufficiently small, $(2\log e)\,k_0(\varepsilon)\lesssim
C_\varepsilon\lesssim \frac{5}{2}k_0(\varepsilon)\log(1/\varepsilon)$.
Consequently, we have $\rho = \frac{1}{2}$, $d_c=2$ and
\beq
\label{37}
L_{\max}(\varepsilon) \lesssim 2^{k_0(\varepsilon)\log(1/\varepsilon)} \leqslant
2^{1/(\pi\sqrt{\varepsilon})\log(1/\varepsilon)}.
\eeq

\begin{remark}
With reference to this last example, the reader interested in sharp bounds on
the $\varepsilon$--capacity in the general setting of Sobolev spaces is referred
to \cite{Birman} (see also Section 6 of \cite{Kolmogorov1}).
\end{remark}

Summarizing, we have the following table:

~

\begin{center}
\begin{tabular}{| c | c | c | c |}
\hline
{\sc Behavior of} $\lambda_k$ & $\log L_{\max}(\varepsilon)$ & $d_c$ & $d_c^{\,e}$ \\
\hline
$\sds e^{\sds -c_1 k}$ & $\sds c'_1 \, [\log(1/\varepsilon)]^2$ & ----- & $2^{1/2}$ \\
\hline
$\sds e^{\sds -c_2 k^2}$ & $\sds c'_2 \, [\log(1/\varepsilon)]^{3/2}$ & ----- & $2^{2/3}$ \\
\hline
$\sds c_3/k^2$ & $\sds c'_3 \, \varepsilon^{-1/2}\log(1/\varepsilon)^{~}$ & $2$ & ----- \\
\hline
\end{tabular}
\vspace{2ex}
\begin{center}
\small
\end{center}
\end{center}

\section{Probabilistic Information}
\label{se:probabilistic}
Here we want to reconsider Equation (\ref{1}) from a probabilistic point of view, adding explicitly the term
representing the noise. With this in mind we pass from Equation (\ref{1}) to Equation (\ref{noisy}), and then to the
probabilistic form of the latter, i.e., Equation (\ref{4}), where $\xi$, $\zeta$ and $\eta$ are Gaussian weak
random variables ($w.r.v.$) in the Hilbert space $L^2[a,b]$ \cite{Balakrishnan}.
A Gaussian $w.r.v.$ is uniquely defined by its mean element and its covariance operator; in the present
case we denote by $R_{\xi \xi}$, $R_{\zeta \zeta}$ and $R_{\eta \eta}$ the covariance operators
of $\xi$, $\zeta$ and $\eta$ respectively. Next, we make the following assumptions:
\begin{itemize}
\item[i)] $\xi$ and $\zeta$ have zero mean, i.e., $m_\xi = m_\zeta = 0$;
\item[ii)] $\xi$ and $\zeta$ are uncorrelated: i.e, $R_{\xi \zeta}$ = 0;
\item[iii)] $R_{\zeta \zeta}^{-1}$ exists.
\end{itemize}
Regarding assumption (i), if it is known that $m_\xi \neq 0$ and $m_\zeta \neq 0$,
then the problem can be easily reformulated
in terms of the variables $(\xi-m_\xi)$ and $(\zeta-m_\zeta)$.
The second hypothesis simply states that the signal--process $\xi$ and the noise--process $\zeta$ are independent.
Finally, the third assumption is the mathematical formulation of the fact that all the components of the data
function are affected by noise or, in other words, that no components of the noise is equal to zero with
probability one. As shown by Franklin, see formula (3.11) of \cite{Franklin}, if assumptions (i) and (ii)
are satisfied, then
\beq
\label{40}
R_{\eta \eta}= A R_{\xi \xi} A^\star + R_{\zeta \zeta},
\eeq
and the cross--covariance operator is given by:
\beq
\label{41}
R_{\xi \eta} = R_{\xi \xi}A^\star.
\eeq
We also assume that $R_{\zeta \zeta}$ depends on a parameter $\varepsilon$ that tends to zero when the noise vanishes, i.e.,
\beq
\label{42}
R_{\zeta \zeta} = \varepsilon^{2} N,
\eeq
where $N$ is a given operator, e.g., $N=I$ for the white noise.

Now, we are faced with the following problem:

\begin{problem}
Given a value $\barg$ of the $w.r.v.$ $\eta$ find an estimate of the $w.r.v.$ $\xi$.
\end{problem}

In order to give an answer to this problem, we turn Equation (\ref{4}) into an infinite sequence of one-dimensional
equations by means of the orthogonal projections, obtaining Equations (\ref{5}), where
$\xi_k = (\xi, \psi_k)$, $\zeta_k = (\zeta, \psi_k)$, $\eta_k = (\eta,\psi_k)$ are
Gaussian random variables. Accordingly we introduce the variances $\rho_k^2 = (R_{\xi \xi} \psi_k, \psi_k)$,
$\varepsilon^2 \nu_k^2 = (R_{\zeta \zeta} \psi_k, \psi_k)$,
$\lambda_k^2 \rho_k^2 + \varepsilon^2 \nu_k^2 = (R_{\eta \eta} \psi_k, \psi_k)$.
Next we evaluate the amount of information on the variable $\xi_k$ which is contained in the variable $\eta_k$;
we have \cite{Gelfand1}:
\beq
\label{44}
J(\xi_k, \eta_k) = - \frac{1}{2} \ln (1 - r_k^2),
\eeq
where
\beq
\label{45}
r_k^2 = \frac{|\EA \left\{ \xi_k \eta_k\right\} |^2}{\EA\left\{ |\xi_k |^2 \right\} \EA\left\{ | \eta_k |^2\right\} } =
\frac{(\lambda_k \rho_k)^2}{(\lambda_k \rho_k )^2 + (\varepsilon \nu_k )^2}.
\eeq
Thus
\beq
\label{46}
J (\xi_k, \eta_k) = \frac{1}{2} \ln \left (1 + \frac{\lambda_k^2 \rho_k^2}{\varepsilon^2 \nu_k^2} \right ).
\eeq
From equality (\ref{46}) it follows that $J (\xi_k, \eta_k) < \frac{1}{2} \ln 2$,
if $\lambda_k\rho_k < \varepsilon \nu_k$, that is if the signal--to--noise ratio of the $k^{\rm th}$ component is small.
Thus, we are naturally led to introduce the following two sets: one, denoted by $\cI$, which accounts for the
components in which the signal dominates the noise; the other one, denoted by $\cN$, which is instead related
to the components in which the noise prevails; precisely, we define:
\begin{eqnarray}
\cI &=& \left\{k\, :\, \lambda_k\rho_k \geqslant \varepsilon\nu_k\right\}, \label{47} \\
\cN &=& \left\{k\, :\, \lambda_k\rho_k < \varepsilon\nu_k\right\}.\label{48}
\end{eqnarray}

\begin{remark}
Let us note that the sets $\cI$ and $\cN$ are not equipped, in general, with any order relation. However, we can
rearrange and renumber the terms $\lambda_k\rho_k$ and $\varepsilon\nu_k$ in such a way as to introduce an order
relationship. Furthermore, for the sake of simplicity and without loss of generality, we hereafter assume that
there do not exist two identical terms $\lambda_k\rho_k/\nu_k$ corresponding to different values of $k$. In this
situation there exists a unique value of $k$, denoted by $k_I$, which separates set $\cI$ from set $\cN$.
\end{remark}

Since $\xi_k$ and $\zeta_k$ are supposed to be Gaussian random variables, we can assume the following probability
densities:
\begin{eqnarray}
p_{\xi_k} (x)& =& \frac{1}{\sqrt{2\pi}\, \rho_k } \exp \left \{ - \left ( \frac{x^2}{2\rho_k^2} \right ) \right\}, \qquad k=1,2, ...,
\label{49} \\
p_{\zeta_k}(x)&=&\frac{1}{\sqrt{2\pi}\,\varepsilon\nu_k}\exp\left\{-\left(\frac{x^2}{2\varepsilon^2\nu_k^2}\right )
\right\}, \qquad k=1,2, .... \label{50}
\end{eqnarray}
By equations (\ref{5}) we can also introduce the conditional probability density $p_{\eta_k}(y|x)$
of the random variable $\eta_k$ for fixed $\xi_k=x$, which reads:
\beq
\label{51}
\begin{split}
p_{\eta_k}(y|x)&=\frac{1}{\sqrt{2\pi}\,\varepsilon\nu_k}\exp\left\{-\frac{(y-\lambda_k x)^2}
{2\varepsilon^2\nu_k^2}\right\} \\
&= \frac{1}{\sqrt{2\pi}\,\varepsilon\nu_k}\exp\left\{-\frac{\lambda_k^2}{2\varepsilon^2\nu_k^2}
\left(x-\frac{y}{\lambda_k}\right )^2\right\}.
\end{split}
\eeq
Let us now apply the Bayes formula, which provides the conditional probability density of $\xi_k$ given
$\eta_k$ through the following expression:
\beq
\label{52}
p_{\xi_k} (x|y) = \frac{p_{\xi_k}(x) p_{\eta_k} (y|x)}{p_{\eta_k} (y)}.
\eeq
Thus, if a realization of the random variable $\eta_k$ is given by $\barg_k$, formula (\ref{52}) becomes
\beq
\label{53}
p_{\xi_k}(x|\barg_k) = A_k \exp\left\{-\frac{x^2}{2\rho_k^2}\right\}
\exp\left\{-\frac{\lambda_k^2}{2\varepsilon^2\nu_k^2}\left(x-\frac{\barg_k}{\lambda_k}\right )^2\right\}, \qquad A_k = \text{\rm Cons.}.
\eeq
The conditional probability density (\ref{53}) can be regarded as the product of two Gaussian probability densities:
\begin{eqnarray*}
p_1(x)&=&A_k^{(1)}\exp\left\{-x^2/2\rho_k^2\right\}, \\
p_2(x) &=& A_k^{(2)}\exp\left\{-(\lambda_k^2/2\varepsilon^2\nu_k^2)\left(x-(\barg_k/\lambda_k)\right)^2\right\},
\end{eqnarray*}
$A_k=A_k^{(1)}\cdot A_k^{(2)}$, whose variances are respectively given by $\rho_k^2$ and $(\varepsilon\nu_k/\lambda_k)^2$.
Let us note that if $k\in\cI$, the variance associated with the density $p_2(x)$ is smaller than the corresponding
variance of $p_1(x)$, and vice versa if $k\in\cN$. Therefore, it is reasonable to consider
as an acceptable approximation of $\langle\xi_k\rangle$ the mean value given by the density $p_2(x)$ if
$k\leqslant k_I$ (i.e., if $k \in \cI$), whereas the mean value given by the density $p_1(x)$ if
$k>k_I$ (i.e., if $k\in\cN$). We can write the following approximation:
\beq
\label{54}
\langle\xi_k\rangle =
\begin{cases}
\sds\frac{\barg_k}{\lambda_k} & \text{if $k \leqslant k_I$}, \\
0 & \text{if $k > k_I$}.
\end{cases}
\eeq
Consequently, given the value $\barg$ of the $w.r.v.$ $\eta$, we are led to consider the following estimate
of $\xi$: $\sum_{k\in\cI}(\barg_k/\lambda_k)\psi_k \equiv \xi_I$. Next, we introduce the operator
$B_\cI:L^2[a,b]\to L^2[a,b]$, defined as follows
\beq
\label{56}
B_\cI \psi_k =
\begin{cases}
\sds\frac{1}{\lambda_k} \psi_k & \text{if $k \leqslant k_I$}, \\
0 & \text{if $k > k_I$},
\end{cases}
\eeq
then $\xi_I=B_\cI\barg=\sum_{k=1}^{k_I}(\barg_k/\lambda_k)\psi_k$. We can now evaluate the global mean
square error; taking into account formulae (\ref{40}) and (\ref{41}), we can formally write:
\beq
\label{58}
\begin{split}
\EA\left\{\|\xi-B_\cI\eta\|^2\right\}&=\Tr\left(R_{\xi\xi}-R_{\xi\xi}AB^\star-BAR_{\xi\xi}+
BR_{\eta\eta}B^\star\right)\\
&=\sum_{k=k_I+1}^\infty \rho_k^2 + \sum_{k=1}^{k_I}\left(\frac{\varepsilon\nu_k}
{\lambda_k}\right)^2.
\end{split}
\eeq
The sum (\ref{58}) is finite if and only if $\Tr R_{\xi\xi}=\sum_{k=1}^\infty\rho_k^2 < \infty$, i.e., if the
covariance operator $R_{\xi\xi}$ is of trace class. In the following we assume that this
condition is satisfied. Hereafter we also suppose that $\lim_{k\to\infty}(\lambda_k\rho_k/\nu_k)=0$,
and therefore the set $\cI$ exists and its cardinality is finite for any given $\varepsilon >0$.
Next, we prove the following lemma.

\begin{lemma}
\label{aux_lemma}
If $\Tr R_{\xi\xi} =\Gamma<\infty$ and moreover $\lim_{k\to\infty}(\lambda_k\rho_k/\nu_k)=0$, then we
can introduce a number $k_\alpha(\varepsilon)$ defined as follows:
\beq
\label{59}
k_\alpha(\varepsilon)=\max\left\{m\in\N\,:\,\sum_{k=1}^m\left(\rho_k^2+\frac{\varepsilon^2\nu_k^2}{\lambda_k^2}
\right)\leqslant\Gamma\right\}.
\eeq
We can then prove:
\begin{eqnarray}
{\rm (i)}&\qquad&\lim_{\varepsilon\to 0} k_\alpha(\varepsilon) = +\infty, \label{60} \\
{\rm (ii)}&\qquad&\lim_{\varepsilon\to 0}\left\{\sum_{k=k_\alpha+1}^\infty\rho_k^2+\sum_{k=1}^{k_\alpha(\varepsilon)}
\left(\frac{\varepsilon\nu_k}{\lambda_k}\right)^2\right\}=0.\label{61}
\end{eqnarray}
\end{lemma}

\begin{proof}
(i) Let us denote by $k_{\alpha_1}$ the sum $k_\alpha+1$. If
equality (\ref{60}) is not true, then there should exist a finite
number $M$, which does not depend on $\varepsilon$ and such that, for
any sequence $\{\varepsilon_i\}$ converging to zero, $k_{\alpha_1} <
M$. From formula (\ref{59}) it then follows:
\beq
\label{62}
\Gamma < \sum_{k=1}^{k_{\alpha_1}(\varepsilon_i)}\left(\rho_k^2+\frac{\varepsilon^2\nu_k^2}{\lambda_k^2}\right)
\leqslant \sum_{k=1}^M
\left(\rho_k^2+\frac{\varepsilon^2\nu_k^2}{\lambda_k^2}\right).
\eeq
For any sequence $\{\varepsilon_i\}$ tending to zero, we have
\beq
\label{63}
\Gamma < \sum_{k=1}^{M} \rho_k^2  \leqslant
\sum_{k=1}^\infty \rho_k^2 = \Gamma,
\eeq
and the contradiction is explicit. \\
(ii) Since $\lim_{\varepsilon\to 0}k_\alpha(\varepsilon)=+\infty$, and $\sum_{k=1}^\infty \rho_k^2 <\infty$, then
\beq
\label{64}
\lim_{\varepsilon\to 0} \sum_{k=k_\alpha(\varepsilon)+1}^\infty \rho_k^2 = 0.
\eeq
Regarding the term $\sum_{k=1}^{k_\alpha(\varepsilon)}(\varepsilon\nu_k/\lambda_k)^2$, we can proceed as follows:
from formula (\ref{59}) we have
\beq
\label{65}
\sum_{k=1}^{k_\alpha(\varepsilon)}\left(\frac{\varepsilon\nu_k}{\lambda_k}\right)^2+
\sum_{k=1}^{k_\alpha(\varepsilon)}\rho_k^2\leqslant\Gamma=\sum_{k=1}^\infty\rho_k^2,
\eeq
and therefore
\beq
\label{66}
\sum_{k=1}^{k_\alpha(\varepsilon)}\left(\frac{\varepsilon\nu_k}{\lambda_k}\right)^2\leqslant
\sum_{k=k_\alpha(\varepsilon)+1}^\infty\rho_k^2.
\eeq
Since $\lim_{\varepsilon\to 0}\sum_{k=k_\alpha+1}^\infty\rho_k^2=0$ (see (\ref{64})), we have
$\lim_{\varepsilon\to 0}\sum_{k=1}^{k_\alpha(\varepsilon)}(\varepsilon\nu_k/\lambda_k)^2=0$.
\end{proof}

Finally, we can prove the following theorem.

\begin{theorem}
\label{theo2}
If the covariance operator $R_{\xi\xi}$ is of trace class, and
$\lim_{\varepsilon\to 0}\lambda_k\rho_k/\nu_k=0$, then the following limit holds true:
\beq
\label{67}
\lim_{\varepsilon\to 0}\EA\left\{\|\xi-B_\cI\eta\|^2\right\}=
\lim_{\varepsilon\to 0}\left\{\sum_{k=k_I+1}^\infty\rho_k^2+\sum_{k=1}^{k_I}\left(\frac{\varepsilon\nu_k}
{\lambda_k}\right)^2\right\}=0.
\eeq
\end{theorem}

\begin{proof}
The proof proceeds in two steps.\\
a) We want to prove that $\lim_{\varepsilon\to 0}\sum_{k=k_I+1}^\infty\rho_k^2=0$.
We have two possibilities: either $k_I\geqslant k_\alpha$, or $k_I < k_\alpha$. In the former case
the statement follows from the fact that $\lim_{\varepsilon\to 0}\sum_{k=k_\alpha+1}^\infty\rho_k^2=0$.
In the latter case, if $k_I < k_\alpha$, then we have:
\beq
\label{68}
\sum_{k=k_I(\varepsilon)+1}^\infty\rho_k^2 \leqslant \sum_{k=k_\alpha+1}^\infty\rho_k^2
+\sum_{k=k_I(\varepsilon)+1}^{k_\alpha(\varepsilon)}\left(\frac{\varepsilon\nu_k}{\lambda_k}\right)^2 \leqslant
\sum_{k=k_\alpha+1}^\infty\rho_k^2 +
\sum_{k=1}^{k_\alpha(\varepsilon)}\left(\frac{\varepsilon\nu_k}{\lambda_k}\right)^2.
\eeq
But in Lemma \ref{aux_lemma} we have proved that the r.h.s. of formula (\ref{68}) tends to zero
as $\varepsilon\to 0$, and the statement follows. \\
b) We want to prove that $\lim_{\varepsilon\to 0}\sum_{k=1}^{k_I(\varepsilon)}\left(\varepsilon\nu_k/\lambda_k\right)^2=0$.
Now again either $k_I \leqslant k_\alpha$ or $k_I > k_\alpha$. In the first case the statement follows from
$\lim_{\varepsilon\to 0}\sum_{k=1}^{k_\alpha(\varepsilon)}\left(\varepsilon\nu_k/\lambda_k\right)^2=0$,
as proved in Lemma \ref{aux_lemma}. If, on the contrary, $k_I > k_\alpha$, then we have, for
$k \leqslant k_I$, $\rho_k \geqslant \varepsilon\nu_k/\lambda_k$, and therefore
\beq
\label{69}
\sum_{k=k_\alpha+1}^{k_I}\left(\frac{\varepsilon\nu_k}{\lambda_k}\right)^2 \leqslant
\sum_{k=k_\alpha+1}^{k_I}\rho_k^2 \leqslant \sum_{k=k_\alpha+1}^\infty\rho_k^2.
\eeq
Since $\lim_{\varepsilon\to 0}\sum_{k=k_\alpha+1}^\infty\rho_k^2=0$, it follows that
\beq
\label{70}
\lim_{\varepsilon\to 0}\sum_{k=k_\alpha+1}^{k_I}\left(\frac{\varepsilon\nu_k}{\lambda_k}\right)^2=0.
\eeq
Now the statement follows recalling that
$\lim_{\varepsilon\to 0}\sum_{k=1}^{k_\alpha}(\varepsilon\nu_k/\lambda_k)^2=0$, as proved in Lemma \ref{aux_lemma}.
\end{proof}

If we now sum up the information carried by the set $\{\eta_k\}_{k\in\cI}$ on the corresponding set
$\{\xi_k\}_{k\in\cI}$ we obtain the quantity:
\beq
\label{71}
\sum_{k=1}^{k_I} J(\xi_k,\eta_k)=
\sum_{k=1}^{k_I}\ln\left(1+\frac{\lambda_k^2\rho_k^2}{\varepsilon^2\nu_k^2}\right)^{1/2} \simeq
\sum_{k=1}^{k_I}\ln\left|\frac{\lambda_k\rho_k}{\varepsilon\nu_k}\right|,
\eeq
which could be called the \emph{probabilistic information} associated with equation (\ref{4}).
For the approximation on the r.h.s. of (\ref{71}) we used $\lambda_k\rho_k\ \geqslant \varepsilon\nu_k$ for $k\in\cI$.
Now, in order to compare the \emph{probabilistic information} with the \emph{metric information},
we may consider two somehow extremal approximations:
\begin{itemize}
\item[$\alpha)$] If $\rho_k \sim \nu_k$, $k\in\cI$, we have
\beq
\label{72}
\sum_{k=1}^{k_I}\ln\left|\frac{\lambda_k\rho_k}{\varepsilon\nu_k}\right| \sim
\sum_{k=1}^{k_I}\ln\left(\frac{\lambda_k}{\varepsilon}\right) =
\sum_{k=1}^{k_0}\ln\left(\frac{\lambda_k}{\varepsilon}\right),
\eeq
since $k_I=k_0$. Let us note that the r.h.s. of formula (\ref{72}) coincides (up to an immaterial conversion factor
between logarithm types) with the lower bound for $H_\varepsilon(\cE)$.
\item[$\beta)$] If $\lambda_k\rho_k \sim \nu_k$, $k\in\cI$, we have
\beq
\label{73}
\sum_{k=1}^{k_I}\ln\left|\frac{\lambda_k\rho_k}{\varepsilon\nu_k}\right| \sim
k_I(\varepsilon)\ln\left(\frac{1}{\varepsilon}\right) \sim k_0(\varepsilon)\ln\left(\frac{1}{\varepsilon}\right),
\eeq
which coincides with the upper bound for $H_{\varepsilon/2}(\cE)$, which we have computed in the various examples of the
previous section.
\end{itemize}

It is interesting to note that the \emph{metric information} provides the limits of the range over which the
probabilistic information varies when the signal--to--noise ratio ranges between the extrema given by the two
previous approximations.  The results given by approximations $(\alpha)$ and $(\beta)$ allow us to look at
the analogy and parallelism between metric and probabilistic information on a more precise and quantitative ground.

\end{document}